\documentclass[10pt]{amsart}   	
									
\usepackage{amsmath}
\usepackage{amssymb, amsfonts}
\usepackage{amsthm}
\usepackage{hyperref}
\usepackage{mathtools}
\usepackage{tikz}
\usepackage{tikz-cd}
\usetikzlibrary{decorations.markings}
\tikzset{help lines/.style={step=#1cm,very thin, color=gray},
help lines/.default=.5} 

\newcommand*\circled[1]{\tikz[baseline=(char.base)]{
            \node[shape=circle,draw,inner sep=2pt] (char) {#1};}}

\usepackage{enumitem}
\usepackage{makecell}
\setcellgapes{4pt}

\usepackage[margin=1in]{geometry}

\newtheorem{thm}{Theorem}[section]
\newtheorem*{thm*}{Theorem}
\newtheorem{lem}[thm]{Lemma}
\newtheorem{cor}[thm]{Corollary}
\newtheorem{prop}[thm]{Proposition}

\newtheorem{innercustomthm}{Theorem}
\newenvironment{customthm}[1]
{\renewcommand\theinnercustomthm{#1}\innercustomthm}{\endinnercustomthm}

\theoremstyle{definition}
\newtheorem{define}[thm]{Definition}

\newtheorem*{define*}{Definition}
\newtheorem{ex}[thm]{Example}

\newtheorem*{nota*}{Notation}
\newtheorem{rem}[thm]{Remark}
\newtheorem*{rem*}{Remark}

\newcommand{\Z}{\mathbb{Z}}

\newcommand{\Hom}{\mathrm{Hom}}
\newcommand{\Ext}{\mathrm{Ext}}
\newcommand{\End}{\mathrm{End}}

\renewcommand{\dim}{\mathrm{dim}}

\renewcommand{\ker}{\mathrm{ker}}

\newcommand{\mods}{\mathsf{mod}}

\usepackage{soul}

\usepackage{tensor}

\title{Resolution quiver and cyclic homology criteria for Nakayama algebras}
\author{Eric J. Hanson, Kiyoshi Igusa}



\begin{document}
\maketitle

\begin{abstract} If a Nakayama algebra is not cyclic, it has finite global dimension. For a cyclic Nakayama algebra, there are many characterizations of when it has finite global dimension. In \cite{shen_homological}, Shen gave such a characterization using Ringel's resolution quiver. In \cite{igusa_cyclic}, the second author, with Zacharia, gave a cyclic homology characterization for when a monomial relation algebra has finite global dimension. We show directly that these criteria are equivalent for all Nakayama algebras. Our comparison result also reproves both characterizations. In a separate paper we discuss an interesting example that came up in our attempt to generalize this comparison result to arbitrary monomial relation algebras \cite{hanson_counterexample}.
\end{abstract} 

\tableofcontents

\section*{Introduction}

A Nakayama algebra is an algebra whose quiver consists of a single oriented cycle with a finite number of monomial relations or a single path, with or without relations. These are called \emph{cyclic} and \emph{linear} Nakayama algebras respectively. There are many results about these algebras. In \cite{ringel_gorenstein}, Ringel introduced the ``resolution quiver'' $R(\Lambda)$ of a Nakayama algebra $\Lambda$ and used it to characterize when $\Lambda$ is Gorenstein and to give a formula for its Gorenstein dimension. In \cite{shen_resolution}, Shen showed that every component of the resolution quiver of a Nakayama algebra has the same weight and, in \cite{shen_homological}, Shen used this to obtain the following result.

\begin{customthm}{A}[Shen \cite{shen_homological}]\label{thm A} A Nakayama algebra has finite global dimension if and only if its resolution quiver has exactly one component and that component has weight 1.
\end{customthm}

It is important for the purpose of induction that this includes linear Nakayama algebras.

In \cite{igusa_cyclic}, the second author and Zacharia showed that a monomial relation algebra has finite global dimension if and only if the relative cyclic homology of its radical is equal to zero. For the case of a cyclic Nakayama algebra, they constructed a finite simplicial complex $L(\Lambda)$, which we will call the ``relation complex'' of $\Lambda$, and show that the reduced homology of $L(\Lambda)$ gives the lowest degree term in the relative cyclic homology of the radical of $\Lambda$. They further simplified their results to obtain the following.

\begin{customthm}{B}[Igusa-Zacharia \cite{igusa_cyclic}]\label{thm B} A Nakayama algebra $\Lambda$ has finite global dimension if and only if the Euler characteristic of the relation complex $L(\Lambda)$ is equal to one.
\end{customthm}

In fact, only cyclic Nakayama algebras were considered in \cite{igusa_cyclic}. However, if $\Lambda$ is a linear Nakayama algebra, it is easy to see that $L(\Lambda)$ is still defined and is contractible. Thus it has Euler characteristic one.

In this paper we will review these definitions and statements and prove the following statement showing that the two results above are equivalent.

\begin{customthm}{C}\label{thm C} For any cyclic Nakayama algebra $\Lambda$ the Euler characteristic of the relation complex $L(\Lambda)$ is equal to the number of components of its resolution quiver $R(\Lambda)$ with weight 1.
\end{customthm}

We prove Theorem \ref{thm C} using an idea that comes from \cite{sen_syzygy}. In fact, our argument proves Theorems \ref{thm A}, \ref{thm B} and \ref{thm C} simultaneously. (See Remark \ref{rem: How the proposition implies Thms ABC}.) We also obtain a new version of Theorem \ref{thm B}:

\begin{customthm}{B'}[Corollary \ref{cor: theorem Bp}]\label{thm Bp} A Nakayama algebra $\Lambda$ has finite global dimension if and only if its relation complex $L(\Lambda)$ is contractible.
\end{customthm}

These results arose in our study, joint with Gordana Todorov, of the concept of ``amalgamation'' \cite{fock_cluster} and the reverse process which we call ``unamalgamation''. Amalgamation is used in \cite{arkani_scattering} to construct plabic diagrams which are used in \cite{alvarez_turaev} to construct new invariants in contact topology.

Since the result of \cite{igusa_cyclic} holds for any monomial relation algebra, we also attempted to generalize the definition of the resolution quiver to such algebras. This lead us to consider examples which eventually led to a counterexample to the $\phi$-dimension conjecture. See \cite{hanson_counterexample} for more details about this story.

\section{Resolution Quiver of a Nakayama Algebra}

By a \emph{Nakayama algebra} of \emph{order} $n$ we mean a finite dimensional algebra $\Lambda$ over a field $K$ given by a quiver with relations where the quiver, which we denote $Q_n$, consists of a single oriented $n$-cycle:
\begin{center}
\begin{tikzcd}
1 \arrow[r,"x_1"] & 2 \arrow[r, "x_2"] & \cdots  \arrow[r, "x_{n-1}"]& n \ar[lll,bend left=30,"x_n" above]
\end{tikzcd}
\end{center}
 with monomial relations given by paths $y_i$ in this quiver written left to right. We write $y_i=x_{k_i}x_{k_i+1}\cdots x_{\ell_i}$ for $i=1,\cdots,r$ where $1\le k_1<k_2<\cdots<k_r\le n$ and $r\ge1$. If the relations have length $\ge2$ they are called \emph{admissible} and the algebra is called a \emph{cyclic Nakayama algebra}.

 We also allow relations of length one. If there is only one of these, say $x_n=0$, we have the linear quiver
 \begin{center}
\begin{tikzcd}
1 \arrow[r,"x_1"] & 2 \arrow[r, "x_2"] & \cdots  \arrow[r, "x_{n-1}"]& n
\end{tikzcd}
\end{center}
modulo admissible relations and the resulting algebra is called a \emph{linear Nakayama algebra}. We also allow more than one relation of length one. If there are $m$ such relations, the algebra is a product of $m$ linear Nakayama algebras arranged in a cyclic order. These are the algebras that we consider under one formalism: $\Lambda=KQ_n/I$ where $KQ_n$ is the path algebra of $Q_n$ and $I$ is the ideal generated by the relations $y_i$. Unless otherwise noted, the relations $y_i$ will be \emph{irredundant}, i.e., no relation $y_i$ is a subword of another relation $y_j$. We sometimes say that the arrows of $Q_n$ are \emph{cyclically composable}.

A finitely generated right $\Lambda$-module $M$ is equivalent to a representation of the quiver $Q_n$ satisfying the relations $y_i=0$. I.e., for each vertex $i$ we have a finite dimensional vector space $M_i=Me_i$ where $e_i$ is the path of length 0 at $i$ and for each arrow $x_i:i\to i+1$ we have a linear map $M_i\to M_{i+1}$ given by the right action of $x_i$.

Note that $\Lambda$ is a graded algebra since all relations are homogeneous. Also, it has been known for a long time that all modules are uniserial. In hindsight we might say this is because $\Lambda$ is a string algebra (see for example \cite{butler_auslander}).

\begin{rem}\label{rem: relations and projectives}
The projective cover $P_j$ of the simple $S_j$ at vertex $j$ is easy to determine given the relations $y_i$: The composition series of $P_j$ starts at $j$ and continues around the quiver $Q_n$ until it completes one of the relations. In particular, the composition series of each projective gives a (possibly redundant) relation for the algebra. The composition series of projectives $P_j$ of length $|P_j|\le |P_{j+1}|$ give the minimal relations.
\end{rem}

\begin{ex}\label{first example} Let $\Lambda_1$ be $KQ_5$ modulo the relations $y_1=x_2x_3$, $y_2=x_3x_4$ and $y_3=x_5x_1x_2$. Then $P_1,P_2$ have composition series ending in $S_2,S_3$; $P_3$ has composition series $S_3,S_4$; and $P_4,P_5$ have composition series ending in $S_5,S_1,S_2$. Composition series are usually written vertically as:
\begin{center}
\begin{tikzpicture}
\node (A) at (0,.8) {$P_1\!:$};
\node at (.5,.8){1};
\node at (.5,.4){2};
\node at (.5,0){3};
\node (B) at (2,.8) {$P_2\!:$};
\node at (2.5,.8){2};
\node at (2.5,.4){3};
\node (C) at (4,.8) {$P_3\!:$};
\node at (4.5,.8){3};
\node at (4.5,.4){4};
\node (D) at (6,.8) {$P_4\!:$};
\node at (6.5,.8){4};
\node at (6.5,.4){5};
\node at (6.5,0){1};
\node at (6.5,-.4){2};
\node (D) at (8,.8) {$P_5\!:$};
\node at (8.5,.8){5};
\node at (8.5,.4){1};
\node at (8.5,0){2};
\end{tikzpicture}
\end{center}
The relations can be read off from the composition series of those projectives $P_i$ so that $|P_i|\le |P_{i+1}|$ which, in this case, are $P_2,P_3,P_5$.
\end{ex}

\begin{define}[Ringel \cite{ringel_gorenstein}]
The \emph{resolution quiver} $R(\Lambda)$ of a Nakayama algebra $\Lambda$ of order $n$ is the quiver having the same vertex set as $Q_n$ and one arrow $i\to j$ if the $i$-th projective module $P_i$ has socle $S_{j-1}$ (so that $j-i$ is the length of $P_i$ modulo $n$). 
\end{define}

\begin{rem}\label{rem: Gustafson's function}
The function $j=f(i)$ used above can be given by
\[
	f(i)=i+|P_i|-n\dim\Hom_\Lambda(P_n,P_i).
\]
This function originates in \cite{gustafson_global} and is called \emph{Gustafson's function}.
\end{rem}

For $\Lambda_1$ from Example \ref{first example}, the lengths of the projectives are $3,2,2,4,3$. So the five arrows of the resolution quiver are $1\to 4, 2\to 4$, $3\to 5, 4\to 3,5\to 3$ giving:
\begin{center}
\begin{tikzcd}
&\circled{1} \ar[rd]\\
R(\Lambda_1):&\circled{2} \arrow[r] & \circled{4}\ar{r} &\circled 3\ar[ r]& \circled 5 \ar[l, bend right=40]
\end{tikzcd}
\end{center}

The \emph{leaves} of the resolution quiver are the vertices which are not targets of arrow. The other vertices, which we call \emph{nodes}, are in bijection with the relations of the algebra.

\begin{define}
The \emph{weight} of an oriented cycle in the resolution quiver is the sum of the lengths of the projective modules $P_i$ for all $i$ in the cycle divided by $n$, the size of the quiver.
\end{define}

For $\Lambda_1$ from Example \ref{first example}, the resolution quiver has one oriented cycle going through vertices $3,5$. The weight is: 
\[
	wt=\frac{|P_3|+|P_5|}n=\frac{2+3}5=1.
\]
Ringel \cite{ringel_gorenstein} showed that each component of the resolution quiver has exactly one oriented cycle and that its weight is a positive integer. The weight of a component of $R(\Lambda)$ is defined to be the weight of its unique oriented cycle. Shen \cite{shen_resolution} proved that all components of $R(\Lambda)$ have the same weight. Then he proved the following.

\begin{thm}[Shen, \cite{shen_homological}]\label{Shen thm}
A cyclic Nakayama algebra $\Lambda$ has finite global dimension if and only if $R(\Lambda)$ has exactly one component with weight one.
\end{thm}

For Example \ref{first example}, $R(\Lambda_1)$ has one component with weight 1. So, $gl\,dim\,\Lambda_1<\infty$.

\begin{ex}\label{second RA example}
Let $\Lambda_2$ be given by the cyclic quiver $Q_5$:

\begin{center}
\begin{tikzcd}
&1 \arrow[r,"x_1"] & 2 \arrow[d, "x_2"] \\
5\ar[ur, bend left, "x_5"] &4\arrow{l}[above]{x_4} & 3 \arrow{l}[above]{x_3}
\end{tikzcd}
\end{center}
modulo the relations: $x_1x_2x_3, x_2x_3x_4x_5, x_4x_5x_1,x_5x_1x_2$ with projectives:

\begin{center}
\begin{tikzpicture}
\node (A) at (0,.8) {$P_1\!:$};
\node at (.5,.8){1};
\node at (.5,.4){2};
\node at (.5,0){3};
\node (B) at (2,.8) {$P_2\!:$};
\node at (2.5,.8){2};
\node at (2.5,.4){3};
\node at (2.5,0){4};
\node at (2.5,-.4){5};
\node (C) at (4,.8) {$P_3\!:$};
\node at (4.5,.8){3};
\node at (4.5,.4){4};
\node at (4.5,0){5};
\node at (4.5,-.4){1};
\node (D) at (6,.8) {$P_4\!:$};
\node at (6.5,.8){4};
\node at (6.5,.4){5};
\node at (6.5,0){1};
\node (D) at (8,.8) {$P_5\!:$};
\node at (8.5,.8){5};
\node at (8.5,.4){1};
\node at (8.5,0){2};
\end{tikzpicture}
\end{center}

Then, in the resolution quiver, we have $1\to 4, 2\to 1, 3\to 2, 4\to 2, 5\to 3$ giving:
\begin{center}
\begin{tikzcd}
R(\Lambda_2): & \circled 5 \ar[r]&\circled{3} \arrow[r] & \circled{2}\ar{r} &\circled 1\ar r & \circled 4 \ar[ll,bend right=30]
\end{tikzcd}
\end{center}

This is connected with weight:
\[
	wt(R(\Lambda_2))=\frac{|P_2|+|P_1|+|P_4|}{n}=\frac{ 4+3+3}5=2
\]
So, by Theorem \ref{Shen thm}, $\Lambda_2$ has infinite global dimension.
\end{ex}

\begin{ex}\label{third RA example}
Let $\Lambda_3$ be $KQ_4$ modulo $rad^2=0$, then all projectives have length 2. So $R(\Lambda_3)$ is:
\begin{center}
\begin{tikzcd}
R(\Lambda_3):&\circled{1} \arrow[r,bend left=25] & \circled{3}\ar[l,bend left=25] &\circled{2} \arrow[r,bend left=25] & \circled{4}\ar[l,bend left=25]
\end{tikzcd}
\end{center}

This has two components, each with weight 1. So again, the algebra has infinite global dimension. In fact, $\Lambda_3$ is self-injective.
\end{ex}

\section{Cyclic Homology of a Monomial Relation Algebra}

When quoting results from \cite{igusa_cyclic}, we note that a cyclic Nakayama algebra is called a ``cycle algebra'' in \cite{igusa_cyclic} and the relation complex is denoted $K$. We changed this to $L$ to avoid confusion with the ground field $K$. We reserve the symbol ``$\Lambda$'' for Nakayama algebras.

For a monomial relation algebra $A$ over any field $K$, the following result was obtained in \cite{igusa_cyclic}.

\begin{thm}\label{thm: characterizing monomial algebras of finite gl dim}
$A$ has finite global dimension if and only if the cyclic homology of the radical of $A$ is zero.
\end{thm}

The statement uses the fact, proved in \cite{igusa_cyclic}, that the relative cyclic homology of the (Jacobson) radical of $A$ (as an ideal in $A$) is equal to the cyclic homology of $J=rad\,A$ as a ring without unit.

We first give the definition of cyclic homology when the ground field $K$ has characteristic zero. (See \cite{igusa_cyclicMatrix} or \cite{loday_cyclic}.) In the proof of Proposition \ref{prop: cyclic complex is good for any char} we outline the more complicated general definition used in \cite{igusa_cyclic}. In both cases we use Proposition 1.2 in \cite{igusa_cyclic}, which allows us to assume that the $X_i$ in the definition below are indecomposable. We use here the term \emph{radical morphism} to mean a morphism $f:X\to Y$ so that $g\circ f$ is nilpotent for any $g:Y\to X$. When $X,Y$ are projective this is equivalent to saying that $f(X)\subseteq rad\,Y$. The radical morphisms $X\to Y$ form a vector space which we denote $r(X,Y)$.

\begin{define}
In the case when the characteristic of the ground field is equal to zero, the cyclic homology $HC_\ast(J)$ of $rad\,A$ is defined as the homology of the \emph{cyclic complex} of $J=rad\,A$ given by
\begin{equation}\label{eq: cyclic complex}
	0\leftarrow C_0(J)\xleftarrow{\overline b} C_1(J)/\Z_2 \xleftarrow{\overline b} C_2(J)/\Z_3 \xleftarrow{\overline b}\cdots
\end{equation}
where $C_p(J)$ is the direct sum of all sequences of $p+1$ cyclically composable radical morphisms between indecomposable projective modules: $X_0\to X_1\to \cdots\to X_p\to X_0$:
\begin{equation}\label{eq: Cp(J)}
	C_p(J)=\coprod_{X_0,\cdots,X_p} r(X_0,X_1)\otimes r(X_1,X_2)\otimes\cdots\otimes r(X_{p-1},X_p)\otimes r(X_p,X_0).
\end{equation}
 $C_p(J)/\Z_{p+1}$ is $C_p(J)$ modulo the action of the cyclic group $\Z_{p+1}$ where the action of the generator $t$ of $\Z_{p+1}$ is given by $t(f_0,f_1,\cdots,f_p)=(-1)^p(f_1,f_2,\cdots,f_p,f_0)$ and $\overline b:C_p(J)/\Z_{p+1}\to C_{p-1}(J)/\Z_p$ is the map induced by $b:C_p(J)\to C_{p-1}(J)$ given by 
\[
b(f_0,\cdots,f_p)=\sum_{i=0}^{p-1} (-1)^{i} (f_0,\cdots,f_if_{i+1},f_{i+2},\cdots,f_p)+(-1)^p(f_pf_0,f_1,\cdots,f_{p-1}).
\]
\end{define}

Since the action of $t\in \Z_{p+1}$ on $C_p(J)$ and the map $b:C_p(J)\to C_{p-1}(J)$ are homogeneous of degree zero, the cyclic homology complex $C_\ast(J)/\Z_\ast$ is graded when $A$ is graded. In the cyclic Nakayama case, we can say more. To form a cycle, the degrees of the morphisms $X_i\to X_{i+1}$ must add up to a multiple of $n$. So $C_\ast(J)/\Z_\ast$ and $HC_\ast(J)$ are graded and have nonzero terms only in degrees which are positive multiples of $n$. Let $HC_\ast^n(J)$ denote the degree $n$ part of $HC_\ast(J)$.

\begin{prop}\label{prop: cyclic complex is good for any char}
The cyclic complex \eqref{eq: cyclic complex} will compute the degree $n$ part of the cyclic homology of $J=rad\,\Lambda$ for $\Lambda$ a Nakayama algebra of order $n$.
\end{prop}

\begin{rem}\label{rem: HC of linear Nakayama is zero}
In the case when $\Lambda$ is a linear Nakayama algebra, the cyclic complex \eqref{eq: cyclic complex} is zero in every degree and the cyclic homology of $J$ is zero in every degree as predicted since linear Nakayama algebras have finite global dimension. 
\end{rem}

\begin{proof}
The cyclic homology of $J$ is in general defined to be the homology of a double complex which is given by replacing $C_p(J)/\Z_{p+1}$ in \eqref{eq: cyclic complex} by a chain complex $C_{p}(J)\otimes_{K\Z_{p+1}} P_\ast(\Z_{p+1})$, where $P_\ast(\Z_{p+1})$ is a free $K\Z_{p+1}$-resolution of $K$. However, in $C_p^n(J)$, the degree $n$ part of $C_p(J)$, the projective modules $X_i$ in \eqref{eq: Cp(J)} must all be distinct. So, the cyclic group $\Z_{p+1}$ acts freely on the set of summands of $C_p^n(J)$. In other words, $C_p^n(J)$ is a free $K\Z_{p+1}$-module. So, the complex $C_{p}^n(J)\otimes_{K\Z_{p+1}} P_\ast(\Z_{p+1})$ has homology only in degree zero where it is $C_p^n(J)/\Z_{p+1}$. Thus the degree $n$ part of the double complex defining $HC_\ast(J)$ collapses to the degree $n$ part of the chain complex \eqref{eq: cyclic complex} as claimed.
\end{proof}

\begin{ex}\label{ex: first calculation of 3rd example}
In Example \ref{third RA example}, $n=4$ and the degree $4$ part of $C_3(J)/Z_4$ contains
\[
	P_4\to P_3\to P_2\to P_1\to P_4
\]
This is a cycle since the composition of any two morphisms is zero. It cannot be a boundary since each map has degree 1. Thus $HC_3(J)\neq0$ showing (again) that $\Lambda_3$ has infinite global dimension.
\end{ex}

In \cite{igusa_cyclic}, a topological interpretation of $HC^n_\ast(J)$ was given. In this definition we use the term ``internal vertex'' of a relation to mean any vertex $i$ so that $x_{i-1}x_i$ is a subword of the relation. For example, the internal vertices of $x_1x_2x_3$ are $2,3$. A relation of length $\ell\le n$ has $\ell-1$ internal vertices. A relation is said to ``cover'' a vertex if that vertex is an interior vertex of the relation.

\begin{define}\label{def: K(A)}
For any (connected) Nakayama algebra $\Lambda$ of order $n$ with $r$ relations we define the \emph{relation complex} to be the simplicial complex $L(\Lambda)$ having one vertex for every relation of length $\le n$. The vertices $v_0,\cdots,v_p$ span a $p$ simplex in $L(\Lambda)$ if and only if the corresponding relations do not cover all $n$ vertices of the quiver.
\end{define}

Let $CL(\Lambda)$ denote the cone of $L(\Lambda)$.

\begin{lem}\cite{igusa_cyclic}\label{lem: HC in terms of L(A)}
For any Nakayama algebra $\Lambda$, the degree $n$ part of the cyclic homology of $rad\,\Lambda$ is equal to the relative homology of $(CL(\Lambda),L(\Lambda))$ with coefficients in $K$:
\[
	HC_p^n(rad\,\Lambda)=H_p(CL(\Lambda),L(\Lambda);K).
\]
\end{lem}

\begin{rem}\label{rem: chi(HC)=1-chi(L)}
This implies in particular that the Euler characteristic of $HC_\ast^n(rad\,\Lambda)$ is $\chi(CL(\Lambda))-\chi(L(\Lambda))=1-\chi(L(\Lambda))$.
\end{rem}

\begin{ex}\label{ex: L(A) is a cone when A is linear}
Suppose $\Lambda$ is a linear Nakayama algebra. In that case one of the relations, say $y_0$, has length 1 and such a relation has no internal vertices. Therefore, $y_0$ can be added to any simplex in $L(\Lambda)$. This means that $L(\Lambda)$ is a cone with cone point $y_0$. So, $L(\Lambda)$ is a contractible space and $H_\ast(CL(\Lambda),L(\Lambda);K)=0$ in every degree. By Remark \ref{rem: HC of linear Nakayama is zero}, this proves Lemma \ref{lem: HC in terms of L(A)} in the linear case.
\end{ex}

\begin{ex}\label{ex: second calculation of 3rd example}
For Example \ref{third RA example}, there are four relations of length 2 giving four vertices $v_0,v_1,v_2,v_3$. Since each relation covers only one vertex, the only set of relations which covers all four vertices is all of them. So, $L(\Lambda_3)$ in this case is a tetrahedron missing its interior. This is homeomorphic to a 2-sphere. So, $HC_p(rad\,\Lambda_3)=H_p(CL(\Lambda_3),K(\Lambda_3))$ is nonzero only when $p=3$, in agreement with the calculation in Example \ref{ex: first calculation of 3rd example}.
\end{ex}

\begin{thm}\cite{igusa_cyclic}\label{thm: HC char of fin gl dim}
For any Nakayama algebra $\Lambda$ the following are equivalent.
\begin{enumerate}
\item $\Lambda$ has finite global dimension.
\item $HC_\ast^n(rad\,\Lambda)=0$.
\item The Euler characteristic of the relation complex $L(\Lambda)$ is 1.
\end{enumerate}
\end{thm}

\begin{ex}\
\begin{enumerate}
\item In Example \ref{first example}, there are two relations neither of which has vertex $4$ in its interior. So, $L(\Lambda_1)$ is two points connected with one edge:
\begin{center}
\begin{tikzpicture}
\coordinate (A) at (0,0);
\coordinate (B) at (1.5,0);
	\draw[fill] (A) circle[radius=2pt];
	\draw[fill] (B) circle[radius=2pt];
	\draw[very thick] (A)--(B);
\end{tikzpicture}
\end{center}
This has Euler characteristic equal to 1. So, $gl\,dim\,\Lambda_1<\infty$.
\item In Example \ref{second RA example}, there are four relations $y_1,y_2,y_3,y_4$ where $y_2=x_2x_3x_4x_5$ and $y_4=x_5x_1x_2$ cover all 5 vertices but the other five pairs do not. Thus $R(\Lambda_2)$ has 4 vertices and 5 edges. One of the two triangles, with vertices $y_1=x_1x_2x_3,y_2,y_3=x_4x_5x_1$ covers $Q_5$. The other ($y_1,y_3,y_4$) does not. So, only the second face is filled in. 
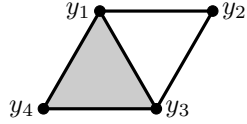
\begin{center}
\begin{tikzpicture}
\coordinate (A) at (0,0);
\coordinate (B) at (1.5,0);
\coordinate (C) at (.75,1.3);
\coordinate (D) at (2.25,1.3);
\draw[fill, color=gray!40!white] (B)--(C)--(A)--cycle;
\foreach \x in {A,B,C,D}
	\draw[fill] (\x) circle[radius=2pt];
	\draw[very thick] (B)--(A)--(C)--(D)--(B)--(C);
	\draw (A) node[left]{$y_4$};
	\draw (C) node[left]{$y_1$};
	\draw (B) node[right]{$y_3$};
	\draw (D) node[right]{$y_2$};
\end{tikzpicture}
\end{center}
So, $\chi L(\Lambda_2)=4-5+1=0$ showing again that $\Lambda_2$ has infinite global dimension.
\item In Example \ref{third RA example}, computed in Example \ref{ex: first calculation of 3rd example} above, $L(\Lambda_3)$ is a 2-sphere with $\chi L(\Lambda_3)=2$. So, again, $gl\,dim\,\Lambda_3=\infty$.
\end{enumerate}
\end{ex}

\section{Comparison of Characterizations}

In this section we prove Theorem \ref{thm C}: For any Nakayama algebra $\Lambda$, the Euler characteristic of its relation complex $L(\Lambda)$ is equal to the number of components of its resolution quiver $R(\Lambda)$ of weight one. In Examples \ref{first example}, \ref{second RA example}, \ref{third RA example} this number is 1,0,2, respectively. The proof will be by induction on the number of vertices of $R(\Lambda)$ which do not lie in any oriented cycle. If this number is zero then $R(\Lambda)$ has no leaves. So, the number of relations is equal to $n$ the size of the quiver $Q_n$. In that case, the relations must all be of the same length, say $\ell$, and we have the following.

\begin{lem}
Let $\Lambda$ be $KQ_n$ modulo $rad^\ell$. Then the number of components of $R(\Lambda)$ is equal to $(n,\ell)$, the greatest common divisor of $n,\ell$, and the weight of each component is $w=\ell/(n,\ell)$. In particular, the weight is 1 if and only if $\ell$ divides $n$ in which case the number of components is $\ell$.
\end{lem}

\begin{proof}
Since the arrows of $R(\Lambda)$ go from $i$ to $i+\ell$ modulo $n$, the size of each oriented cycle is the smallest number $m$ so that $m\ell$ is a multiple of $n$. This number is $m=n/(n,\ell)$. The weight of the cycle is
\[
	wt=\frac{m\ell}n=\frac{\ell}{(n,\ell)}.
\]
Since the total length of all the projectives is $n\ell$, the number of components is $n\ell/m\ell=n/m=(n,\ell)$ as claimed.
\end{proof}

\begin{lem}
Let $\Lambda$ be $KQ_n$ modulo $rad^\ell$. Then the Euler characteristic of $L(\Lambda)$ is
\[
	\chi(L(\Lambda))=\begin{cases} \ell & \text{if $\ell$ divides $n$} \\
   0 & \text{otherwise}
    \end{cases}
\]
\end{lem}

\begin{proof}
By Remark \ref{rem: chi(HC)=1-chi(L)}, this is equivalent to the equation proved in Theorem 4.4 of \cite{igusa_cyclic} that $\chi(HC_\ast^n(rad\,\Lambda))=1-\ell$ if $\ell$ divides $n$ and is 1 otherwise. This includes the degenerate case when $\ell=1$ since, in this case, $L(\Lambda)$ is contractible by Example \ref{ex: L(A) is a cone when A is linear} and thus $\chi(L(\Lambda))=1=\ell$.
\end{proof}

These two lemmas imply that Theorem \ref{thm C} holds in the case when $R(\Lambda)$ is a union of oriented cycles or, equivalently, when $R(\Lambda)$ has no leaves. These lemmas also imply Theorems \ref{thm A} and \ref{thm B} in this case since $KQ_n$ modulo $rad^\ell$ has finite global dimension if and only if $\ell=1$.

Suppose now that $\Lambda$ is a Nakayama algebra of order $n$ and $j$ is a leaf of $R(\Lambda)$. By reindexing we may assume $j=n$.

\begin{prop}\label{prop: properties of L'} Let $n$ be a leaf of $R(\Lambda)$ and let $\Lambda'=\End_\Lambda(P')$, where $P'=P_1\oplus \cdots \oplus P_{n-1}$. Then we have the following.
\begin{enumerate}
\item $R(\Lambda')$ is equal to $R(\Lambda)$ with the leaf $n$ removed. In particular:
	\begin{enumerate} 
	\item $R(\Lambda')$ has the same number of components as $R(\Lambda)$ and
	\item $R(\Lambda')$ has one fewer vertex not in an oriented cycle than $R(\Lambda)$.
	\end{enumerate}
\item $R(\Lambda')$ has the same weight as $R(\Lambda)$.
\item $L(\Lambda')$ is homotopy equivalent to $L(\Lambda)$. So, they have the same Euler characteristic.
\item $gl\,dim\,\Lambda'\le gl\,dim\,\Lambda\le gl\,dim\,\Lambda'+2 $. Thus, $gl\,dim\,\Lambda<\infty$ if and only if $gl\,dim\,\Lambda'<\infty$.
\end{enumerate}
\end{prop}

\begin{rem}\label{rem: How the proposition implies Thms ABC}
Suppose for a moment that $\Lambda'$ has all of these properties. The proof of Theorem \ref{thm C} then proceeds as follows. By induction on the number of vertices of the resolution quiver which are not in any oriented cycle, Theorem \ref{thm C} holds for $\Lambda'$. So $\chi(L(\Lambda'))$ is equal to the number of components of $R(\Lambda')$ with weight 1. But these numbers are the same for $\Lambda'$ and $\Lambda$. So, Theorem \ref{thm C} also holds for $\Lambda$. By (4), Theorems \ref{thm A}, \ref{thm B} also hold for $\Lambda$ assuming they hold for $\Lambda'$. Therefore, the existence of $\Lambda'$ satisfying the properties above will prove all three theorems at the same time!
\end{rem}

\subsection{Algebraic properties of $\Lambda'$}

Since $P'$ is a projective $\Lambda$-module, we have an exact functor
\[
	\pi: \mods\text-\Lambda\to \mods\text-\Lambda'
\]
given by $\pi(M_\Lambda)=\Hom_\Lambda(\prescript{}{\Lambda'}{P'_\Lambda},M_\Lambda)$.

\begin{lem}\label{lem: id Sn=1}
The injective dimension of $S_n$ is at most one.
\end{lem}

\begin{proof}
Since $n$ is a leaf in $R(\Lambda)$, there is no relation for $\Lambda$ which ends in $x_{n-1}$. So, the injective envelope of $S_n$ modulo $S_n$ is injective making $S_n\to I_n\to I_n/S_n\to 0$ the injective copresentation of $S_n$.
\end{proof}

The following result appears in \cite{ChenYe_Gorenstein} and is one step in the construction of $\varepsilon(\Lambda)$ in \cite{sen_syzygy}.

\begin{prop}\label{prop: embedding of Lambda-prime}
The functor 
\[
	\alpha:\mods\text-\Lambda'\to \mods\text-\Lambda
\]
given by $\alpha(X_{\Lambda'})=X\otimes_{\Lambda'}P'$ is a full and faithful exact embedding which takes projectives to projectives and whose image contains the second syzygy of any $\Lambda$-module.
\end{prop}

\begin{proof} Since $\alpha(\Lambda')=P'$ and $\alpha$ is additive, $\alpha$ takes projective modules to projective modules. 

For any $X\in \mods\text-\Lambda'$, consider a free presentation of $X$ of the form 
\begin{equation}\label{eq: presentation of X}
	(\Lambda')^p\xrightarrow f (\Lambda')^q\to X\to 0.
\end{equation}
 Since $\alpha$ is right exact, we get a projective presentation of $\alpha(X)$: 
\begin{equation}\label{eq: presentation of aX}
	(P')^p\xrightarrow{f_\ast} (P')^q\to \alpha(X)\to 0.
\end{equation}
Since $\Hom_\Lambda(P',P')=\Hom_{\Lambda'}(\Lambda',\Lambda')$, the correspondence $f\leftrightarrow f_\ast$ is a bijection with $\pi$ sending \eqref{eq: presentation of aX} back to \eqref{eq: presentation of X}. If we do the same for another $\Lambda'$-module $Y$, the morphisms between $P'$-presentations in $\mods\text-\Lambda$ are in bijection with the morphisms between free presentations in $\mods\text-\Lambda'$. So, we get \[
\Hom_\Lambda(\alpha(X),\alpha(Y))=\Hom_{\Lambda'}(X,Y)\]
showing that $\alpha$ is full and faithful.

The image of the functor $\alpha$ contains all cokernels of all morphisms $(P')^p\to (P')^q$ as in \eqref{eq: presentation of aX}. Since $id\,S_n\le 1$, this includes all second syzygies of all $\Lambda$-modules.

To show that $\alpha$ is exact, consider a short exact sequence $0\to Z\to Y\to X\to 0$ in $\mods\text-\Lambda'$. Take free presentations of $X,Z$, and combine them to get a presentation of $Y$:
\begin{center}
\begin{tikzcd}
0 \ar[r] & (\Lambda')^{r}\ar[r]\ar[d,"h"] & (\Lambda')^{r+p}\ar[r]\ar[d,"g"]& (\Lambda')^{p}\ar[r]\ar[d,"f"] &0\\
0 \ar[r] & (\Lambda')^{s}\ar[r]\ar[d] & (\Lambda')^{s+q}\ar[r]\ar[d]& (\Lambda')^{q}\ar[r]\ar[d] &0\\
0 \ar[r]  & Z \ar[r]\ar[d] & Y \ar[r]\ar[d] & X\ar[r]\ar[d] & 0\\
&0&0&0
\end{tikzcd}
\end{center}
Applying the functor $\alpha$ we get a right exact sequence:
\[
	\alpha(Z)\xrightarrow k \alpha(Y)\to \alpha(X)\to 0.
\]
Since the exact functor $\pi=\Hom_\Lambda(P',-)$ sends this back to the sequence $0\to Z\to Y\to X\to 0$, $\ker\,k$ must be a direct sum of copies of $S_n$. However, by the Snake Lemma, the kernel of $f_\ast:(P')^p\to (P')^q$ maps onto $\ker\,k$. But this is impossible unless $\ker\,k = 0$, since $\Ext^2_\Lambda(\alpha(X),S_n)=0$ by Lemma \ref{lem: id Sn=1}. Therefore, $\alpha$ is exact as claimed.
\end{proof}

\begin{cor}\label{cor: gl dim of Lambda and Lambda-prime}
The global dimensions of $\Lambda$ and $\Lambda'$ differ by at most two:
\[
gl\,dim\,\Lambda'\le gl\,dim\,\Lambda\le gl\,dim\,\Lambda'+2.\]
In particular, $gl\,dim\,\Lambda<\infty$ if and only if $gl\,dim\,\Lambda'<\infty$.
\end{cor}

\begin{proof} Since $\alpha$ is an exact functor taking projectives to projectives, it takes a projective resolution of $M$ to a projective resolution of $\alpha(M)$. Since $\alpha$ is full and faithful, the projective resolution of $M$ does not split. This means $M, \alpha(M)$ have the same projective dimension. So,
\[
	gl\,dim\, \Lambda'\le gl\,dim\,\Lambda.
\]
Given any $\Lambda$-module $M$ of projective dimension at least 3, by Proposition \ref{prop: embedding of Lambda-prime} there is a $\Lambda'$-module $N$ so that $\alpha(N)\equiv \Omega^2M$. Thus, \[
pd\,M=2+pd\,N\le gl\,dim\,\Lambda'\]
which implies that $gl\,dim\,\Lambda\le 2+gl\,dim\,\Lambda'$.
\end{proof}

The following essentially follows from Lemma 2.1 in \cite{shen_homological}.

\begin{prop}
The resolution quiver of $\Lambda'$ is equal to the resolution quiver of $\Lambda$ with the leaf $n$ removed. In particular, $R(\Lambda),R(\Lambda')$ have the same number of components. Furthermore, the weight of $R(\Lambda')$ is equal to the weight of $R(\Lambda)$.
\end{prop}

\begin{proof}
The length of $\pi(M)$ for any $\Lambda$-module $M$ is $|M|-\dim\Hom_\Lambda(P_n,M)$. Gustafson's function for $\Lambda$ on any $i<n$ is
\[
	f_{\Lambda}(i)=i+|P_i|-n\dim \Hom_\Lambda(P_n,P_i)=i+|\pi(P_i)|-(n-1)\Hom_\Lambda(P_n,P_i),
\]
which is congruent to $f_{\Lambda'}(i)$ modulo $n-1$. Thus, for $i<n$, $i\to j$ in $R(\Lambda)$ if and only if $i\to j$ in $R(\Lambda')$. In particular, the cycles are the same.

Suppose that $j_1\to j_2\to \cdots \to j_m$ is a cycle in $R(\Lambda)$ with weight
\[
	w=\frac1n \sum |P_{j_i}|.
\]
This implies that $\bigoplus P_{j_i}$ has each simple $\Lambda$-module $w$ times in its composition series. In particular, $w=\sum \dim\Hom(P_n,P_{j_i})$. The same cycle is also a cycle in $R(\Lambda')$ with weight
\[
	w=\frac1{n-1} \sum |\alpha(P_{j_i})|=\frac1{n-1} \sum \left(
	|P_{j_i}| -\dim\Hom_\Lambda(P_n,P_{j_i})
	\right)= \frac{nw-w}{n-1}=w.
\]
So, $R(\Lambda)$ and $R(\Lambda')$ have the same weight.
\end{proof}

\subsection{Topology of $L(\Lambda')$} It remains to show that the relation complex of $\Lambda'$ is homotopy equivalent to that of $\Lambda$. This will follow from the following three lemmas where we note that only Lemma \ref{lem: redundant relations give the same L(L)} requires $n$ to be a leaf of the resolution quiver.

For any word $y$ in the letters $x_1,\cdots,x_n$ let $\delta(y)$ be $y$ with all occurrences of the letter $x_n$ deleted.

\begin{lem}\label{lem: relations for L'} Given that $\Lambda=KQ_n$ modulo relations $y_1,\cdots,y_r$ with $y_i=x_{j_i}x_{j_i+1}\cdots x_{\ell_i}$, $\Lambda'$ is equal to $KQ_{n-1}$ modulo the possibly redundant relations $y_i'$ given by
\[
	y_i'=\begin{cases} x_{n-1}\delta(y_i) & \text{if }j_i=n\\
    \delta(y_i)& \text{otherwise}
    \end{cases}
\]
\end{lem}

\begin{proof}
We use the fact that the composition series of any projective module is a relation for the algebra and that the minimal relations come from projective modules $P_i$ with length $|P_i|\le |P_{i+1}|$.

If $j_i\neq n$ then $y_i$ gives the composition series of $P_{j_i}$ and $y_i'=\delta(y_i)$ gives the composition series of $\pi(P_{j_i})$. If $j_i=n$ then $y_i$ gives the composition series for $P_n$. Since $P_n$ maps onto the radical of $P_{n-1}$, $x_{n-1}y_i$ is a relation for $\Lambda$ and $y_i'=x_{n-1}\delta(y_i)$ is a relation for $\Lambda'$.

We conclude that $y_1',\cdots,y_r'$ are among the relations that hold for $\Lambda'$. We need to show that this list contains all the minimal relations for $\Lambda'$. So, suppose that $|\pi(P_i)|\le |\pi(P_{i+1})|$ and $z$ is the relation for $\Lambda'$ given by $\pi(P_i)$. If $i\le n-2$ then $|P_i|\le |P_{i+1}|$; so, $z=\delta(y)=y'$, where $y$ is the relation for $\Lambda$ given by $P_i$. Thus, $z$ is in the list. If $i=n-1$ then $|\pi(P_{n-1})|\le |\pi(P_1)|$, which implies that $|P_{n-1}|\le |P_1|+1$. So, either
\begin{enumerate}
\item[(a)] $|P_{n-1}|\le |P_n|$ in which case $z=\delta(y)=y'$, where $y$ is the relation for $\Lambda$ given by $P_{n-1}$ or
\item[(b)] $|P_n|=|P_{n-1}|-1\le |P_1|$, so $P_n$ gives a relation $y$ for $\Lambda$ and $y'=x_{n-1}\delta(y)=z$.
\end{enumerate}
In both cases the relation $z$ is in the list. So, $\{y_i'\}$ is a list of relations for $\Lambda'$ which contains all the minimal relations.
\end{proof}

\begin{lem}\label{lem: redundant relations give the same L(L)}
Let $L'(\Lambda')$ be the simplicial complex with vertices given by the possibly redundant relations $y_i'$ of Lemma \ref{lem: relations for L'} with $p+1$ such relations spanning a $p$-simplex if their interiors cover all of $Q_{n-1}$. Then $L'(\Lambda')$ is isomorphic to $L(\Lambda)$ as a simplicial complex.
\end{lem}

\begin{proof}
We have a bijection $y_i\leftrightarrow y_i'$ between the vertices of $L(\Lambda)$ and those of $L(\Lambda')$. So, we need to show that a collection of relations, $\{y_i\}$ for $\Lambda$ covers all the vertices of $Q_n$ if and only if the corresponding relations $y'_i$ cover all the vertices of $Q_{n-1}$. We claim the following.
\begin{enumerate}
\item If $y'_i$ covers vertex $n-1$ then $y_i$ covers $n-1$ and $n$. 
\item $y_i$ covers vertex $k\neq n$ if and only if $y_i'$ covers $k$.
\end{enumerate}
To prove (1), suppose that $y'_i$ covers $n-1$. Then this relation contains $x_{n-2}x_{n-1}$ as a subword. Since $n$ is a leaf of $R(\Lambda)$, the relation $y_i$ cannot end in $x_{n-1}$. Thus, it must contain the subword $x_{n-2}x_{n-1}x_n$. So, $y_i$ covers $n-1$ and $n$ as claimed.

To prove (2), first suppose $y_i$ covers $k\neq n$. Then $\delta(y_i)$ will cover $k$ except in the case when $k=1$ and $x_nx_1$ are the first two letters of $y_i$. In that case $y_i'=x_{n-1}x_1\cdots$ still covers $k=1$. The converse follows from this and (1). This proves the lemma.
\end{proof}

\begin{lem}\label{lem: removing redundant relations is collapse}
The relation complex $L(\Lambda')$ is homotopy equivalent to one constructed using any redundant set of relations such as $L'(\Lambda')$.
\end{lem}

\begin{proof}
Suppose that the relation set used to construct $L'(\Lambda')$ has relations $z_0$ and $z$ so that $z_0$ is a subword of $z$. Let $L''(\Lambda')$ be obtained from $L'(\Lambda')$ by deleting the vertex $z$ and all open simplices containing $z$ as an endpoint. The union of these is called the ``open star'' of $z$ and the closure of this union is called the ``closed star'', $St(z)$ (see \cite[Chapter 2.C]{hatcher_algebraic}). The union of all simplices in $St(z)$ which do not contain $z$ as an endpoint is called the ``link'' of $z$ and denoted $Lk(z)$. Then
\[
	L'(\Lambda')=L''(\Lambda')\cup_{Lk(z)}St(z),
\]
which means that $L'(\Lambda')=L''(\Lambda')\cup St(z)$ and $L''(\Lambda')\cap St(z)=Lk(z)$. Since $St(z)$ is a cone on $Lk(z)$, it is contractible. So, $L'(\Lambda')\simeq L''(\Lambda')$ if and only if $Lk(z)$ is contractible.

But, $St(z)$ is the cone on $z_0$ in $L'(\Lambda')$ since, for any set of relations including $z$ which does not cover all vertices of the quiver, adding $z_0$ will cover the same set of vertices (because $z_0$ is a subword of $z$). This implies that the link of $z$ is the cone on $z_0$ in $L''(\Lambda')$. So, $L'(\Lambda')\simeq L''(\Lambda')$. Repeat this to remove redundant relations one at a time without changing the homotopy type of $L'(\Lambda')$.
\end{proof}

\begin{ex}
For $\Lambda_1$ in Example \ref{first example}, $n=5$ is the only leaf. The relations $y_i$, $y_i'$ are given by:
\[
\begin{array}{c|cc}
y_i & y_i'\\
\hline
x_1x_2x_3 & x_1x_2x_3\\
x_2x_3x_4x_5 & x_2x_3x_4\\
x_4x_5x_1 & x_4x_1 & \text{($z_0$: subword of $z$)}\\
x_5x_1x_2 & x_4x_1x_2 & \text{($z$: redundant)}
\end{array}
\]
The relation complex of $\Lambda_1'$ is obtained from $L'(\Lambda_1')=L(\Lambda_1)$ by deleting all open simplices containing $z=y_4'=x_4x_1x_2$ as a vertex.  So, $L(\Lambda_1')$ is the hollow triangle with vertices $y_1',y_2',y_3'$. One can see that this is a deformation retract of $L(\Lambda_1)$. The proof of Lemma \ref{lem: removing redundant relations is collapse} uses the fact that $Lk(z)=Lk(y_4')$ is a cone on $z_0=y'_3$ and is thus contractible. So, $L(\Lambda_1')$ is homotopy equivalent to $L(\Lambda_1)$.
\begin{center}
\begin{tikzpicture}

\begin{scope}
\coordinate (X) at (-2,.7);
\draw (X) node{$L'(\Lambda_1')=L(\Lambda_1):$};
\coordinate (A) at (0,0);
\coordinate (B) at (1.5,0);
\coordinate (C) at (.75,1.3);
\coordinate (D) at (2.25,1.3);
\coordinate (BC) at (1.2,.6);
\coordinate (BC1) at (2,1);
\coordinate (BC2) at (.9,0);
\coordinate (BC3) at (2.5,.5);
\draw[thick, color=blue,<-] (BC)..controls (BC1) and (BC2)..(BC3);
\draw[color=blue] (BC3) node[right]{$Lk(z)$}; 
\draw[fill, color=gray!40!white] (B)--(C)--(A)--cycle;
\foreach \x in {A,D}
	\draw[fill] (\x) circle[radius=2pt];
	\draw[very thick] (B)--(A)--(C)--(D)--(B);
	\draw[very thick, color=blue] (B)--(C);
\foreach \x in {B,C}
	\draw[fill,color=blue] (\x) circle[radius=2pt];
	\draw (A) node[left]{$z=y_4'$};
	\draw (C) node[left]{$y'_1$};
	\draw (B) node[right]{$z_0=y'_3$};
	\draw (D) node[right]{$y'_2$};
\end{scope}
\begin{scope}[xshift=6cm]
\coordinate (X) at (-.6,.7);
\draw (X) node{$L(\Lambda_1'):$};
\coordinate (A) at (0,0);
\coordinate (B) at (1.5,0);
\coordinate (C) at (.75,1.3);
\coordinate (D) at (2.25,1.3);
\foreach \x in {C,B,D}
	\draw[fill] (\x) circle[radius=2pt];
	\draw[very thick] (B)--(D)--(C)--cycle;
	\draw (C) node[left]{$y'_1$};
	\draw (B) node[right]{$z_0=y'_3$};
	\draw (D) node[right]{$y'_2$};
\end{scope}
\end{tikzpicture}
\end{center}
\end{ex}

From Lemmas \ref{lem: redundant relations give the same L(L)} and \ref{lem: removing redundant relations is collapse} we immediately obtain the following.

\begin{prop}\label{prop: L simeq Lp}
The relation complexes $L(\Lambda)$ and $L(\Lambda')$ are homotopy equivalent.\qed
\end{prop}

This concludes the proof of Proposition \ref{prop: properties of L'} which proves Theorems \ref{thm A}, \ref{thm B} and \ref{thm C} simultaneously. (See Remark \ref{rem: How the proposition implies Thms ABC}.) Theorem \ref{thm Bp} also follows:

\begin{cor}\label{cor: theorem Bp}
A Nakayama algebra has finite global dimension if and only if its relation complex is contractible.
\end{cor}

\begin{proof}
If $L(\Lambda)$ is contractible, its Euler characteristic is 1 and it thus has finite global dimension by Theorem \ref{thm C}. Conversely, if $\Lambda$ has finite global dimension then, by repeatedly removing the leaves of the resolution quiver we obtain an algebra $\Lambda''$ with $rad=0$, i.e., it is semi-simple. Since the relations of $\Lambda''$ have empty interiors, the relation complex of $\Lambda''$ is a solid simplex, so it is contractible. By Proposition \ref{prop: L simeq Lp}, $L(\Lambda)\simeq L(\Lambda'')$ is also contractible.
\end{proof}

\section{Comments}

The construction of $\Lambda'$ from $\Lambda$ is an example of ``unamalgamation'', which is the construction described in Lemma \ref{lem: relations for L'}. This is a partial inverse to the ``amagamation'' construction of Fock-Goncharov \cite{fock_cluster} used in \cite{arkani_scattering} to construct the Jacobian algebras of plabic diagrams \cite{postnikov_total}. These operations are partial inverses since we claim that $UAU=U$ and $AUA=A$, where $A,U$ stand for amalgamation and unamalgamation, respectively.

\begin{center}
\begin{tikzcd}
& 1' \ar[d]&1 \arrow[r,"x_1"] & 2 \arrow[d, "x_2"] \\
5\ar[ur, bend left, "x_5"] &4'\arrow{l}[above]{x_4} &4\arrow{u} & 3 \arrow{l}[above]{x_3}
\end{tikzcd} 
$\text{amalgamation }\Rightarrow\over
\Leftarrow\text{ unamalgamation}$
\begin{tikzcd}
&1 \arrow[r,"x_1"] & 2 \arrow[d, "x_2"] \\
5\ar[ur, bend left, "x_5"] &4\arrow{l}[above]{x_4} & 3 \arrow{l}[above]{x_3}
\end{tikzcd} 
\end{center}

The idea to consider removing the leaves $R(\Lambda)$ comes from \cite{sen_syzygy} where, for $\Lambda$ of infinite projective dimension, Emre \c{S}en constructs a Nakayama algebra $\varepsilon(\Lambda)$ which seems to have the property that $R(\varepsilon(\Lambda))$ is equal to $R(\Lambda)$ with all its leaves removed. He then proves that the $\phi$-dimension of $\Lambda$ is equal to 2 plus the $\phi$-dimension of $\varepsilon(\Lambda)$.

\section*{Acknowledgements}
Both authors thank Gordana Todorov who is working with them on the larger project of amalgamation and unamalgamation of Nakayama algebras. The first author thanks Job Rock for meaningful conversations. The second author thanks Emre \c{S}en for many discussions about Nakayama algebras and especially the construction of $\varepsilon(\Lambda)$ from \cite{sen_syzygy}. The second author thanks his other coauthor Daniel \'Alvarez-Gavela for working with him, in \cite{alvarez_turaev}, to reduce the calculation of higher Reidemeister torsion in \cite{IK93} to a calculations on plabic diagrams (\cite{postnikov_total}). Also, the second author would like to thank An Huang for directing him to the paper \cite{arkani_scattering} which inspired this line of research, particularly the transition from \cite{alvarez_turaev} to \cite{hanson_counterexample} and this paper. Finally, the authors thank the referee for helpful comments.

\newcommand{\etalchar}[1]{$^{#1}$}

\end{document}